\documentclass[12pt]{amsart}
\usepackage{graphicx}
\newtheorem{thm}{Theorem}[section]

\newtheorem{lem}[thm]{Lemma}

\theoremstyle{definition}

\theoremstyle{remark}

\numberwithin{equation}{section}

\begin{document}

\title[Fubini-Tonelli type theorem]{Fubini-Tonelli type theorem for non product measures  in a product space}%
\author{ Jorge Salazar}

\address{DMAT, Universidade de \'Evora, \'Evora - Portugal}%
\email{salazar@uevora.pt}%

\subjclass{ }%

\date{}

\begin{abstract}
I prove a theorem about iterated integrals for non-product measures in a product space. The first task is to show the existence of a family of measures on the second space, indexed by the points on of the first space (outside a negligible set), such that integrating the measures on the index against the first marginal gives back the original measure (see Theorem \ref{teo}). At the end, I give a simple application in Optimal Transport.
\end{abstract}
\maketitle
\section{Introduction}
\noindent
The Fubini-Tonelli theorem states that the integral of a function defined on a product space, against a measure which is a product of  measures on the factor spaces, can be obtained  by iterated  integration, i.e. integrating one variable  (against its marginal measure) at the time. 
\vskip 2mm
\noindent
 If a measure on a product space is not a product measure, is it still possible to decompose the measure and evaluate the integral using iterated integration? To better understand the problem, imaging we are dealing with a measure $ \zeta  $ which is absolutely continuous with respect to the product measure $ \mu\otimes\nu $, i.e. there is a function $ \delta: X\times Y\rightarrow \left[ 0, \infty\right] $, such that for all measurable set $C\subseteq X\times Y$,
 \[ \zeta\left( C\right) =\int_{C} \delta(x,y)\,\mu\otimes\nu\left( dx, dy \right) . \]
 \vskip 2mm
 \noindent
 Then, using the classical Fubini-Tonelli theorem, we can decompose this integral into 
 \[ \zeta\left( C\right) =\int_X\left( \int_{C_x} \delta(x,y)\, \nu\left( dy \right) \right)  \mu\left( dx \right) ,\]
 where  $C_x :=\left\lbrace y\in Y;\, \left( x,y\right)\in C \right\rbrace  $ 
  is the slice of $C$ at $X$.
 Accordingly, $\nu$ decomposes into the measures
 \[ \nu_x(dy):=\delta(x,y)\, \nu\left( dy \right) , \]
 which integrates against $\mu$ to give  $\zeta$. Symbolically,
 \[ \zeta\left( dx, dy \right) =\nu_x(dy)\, \mu\left( dx \right) . \]
  With this decomposition the order of integration is not interchangeable, since the first measure depends on the second variable.
 To interchange the order of integration, we must decompose  $\mu$ in a similar way,
  \[ \mu_y(dx):=\delta(x,y)\, \mu\left( dx \right) , \]
  which integrates against $\nu$ to give  $\zeta$. Symbolically,
  \[ \zeta\left( dx, dy \right) =\mu_y(dx)\, \nu\left( dy \right) . \]
\noindent
 In this paper,  the existence of this kind of decomposition is established for arbitrary  Borel probability measures on 
 the product of two complete, separable,  locally compact, 
 metric spaces  (see Theorem \ref{teo}). I restricted myself  to the case of probability measures to simplify the discourse, although the results stay valid for $\sigma-$finite measures.
\vskip 2mm
\noindent
 In the best of my knowledge, there is nothing of the kind in the literature on foundations of measure theory that describes similar results. Nonetheless, this question is natural and I believe it may provide a useful calculation and/or analytical tool as much as the classical Fubini-Tonelli theorem does. 
 \vskip 2mm
 \noindent
 Optimal Transport, for example,  deals with fixed marginal probability measures and a minimal cost is seek  among all the couplings of the given marginal probabilities, i.e. among all the probabilities on the product space, such that the marginal measures are the ones given.
 It would be a nice research project to look for a new characterization of the optimal transport plans in terms of the measures along the ``fibers'', obtained from  the  decomposition described in Theorem \ref{teo}.
In section \ref{OT}, I give a simple application of Theorem \ref{teo}, showing that a pair of competitive price functions, whose integral with respect to some transference plan matches the transport cost, are conjugate to each other almost surely. This complements the Kantorovich duality theorem on the nature, regarding convexity/concavity, of pair of competitive prices maximizing the profit. See Villani's book \cite{vill}, page 70, for a very detailed discussion of the Kantorovich theorem. 
  \vskip 2mm
  \noindent
Another interesting project is the application of Theorem \ref{teo} to the study of measures on the Tangent bundle of  Riemannian manifolds. Indeed, the local charts of the tangent bundle are Cartesian products of Euclidean  open sets. Using local charts,   we can transport the measure to this product to be decomposed and then sent back the family of measures fiber-wise. In the literature, the measures on tangent bundles are a kind of product measures, as is the volume  obtained from the  Sasaki \cite{sasaki}  metric, or the measure on the unit sphere on the tangent space, integrated against the volume element of the base manifold. I believe Theorem \ref{teo}  is a tool that could help  exploring general integration on tangent bundles.

\section{Main theorem}

{\thm\label{teo} Let $X\times Y$ be the product of two complete, separable,  locally compact, 
	metric spaces. We equip $X$, $Y$, and $X\times Y$  with their Borel $\sigma-$algebras, denoted by $\mathfrak{B}_X$, $\mathfrak{B}_Y$, and $\mathfrak{B}_{X\times Y}$ respectively.  
	\vskip 3pt
	\noindent 
	Let $\zeta$ be a  probability measure on $\mathfrak{B}_{X\times Y}$, and denote $\mu$ and $\nu$  the marginal probabilities on $\mathfrak{B}_X$, $\mathfrak{B}_Y$ respectively. i.e.
	\[ \forall A\in \mathfrak{B}_X, \ \mu(A)= \zeta(A\times Y) \]
	and 
	\[ \forall B\in \mathfrak{B}_Y, \ \nu(B)= \zeta(X\times B) .\]
\vskip 3pt
\noindent 
Then, outside an exceptional $\mu-$negligible set $E_1 \in \mathfrak{B}_X$  ($\mu\left( E_1\right) = 0$), for all 
	$ x \in X\setminus E_1 $, there is a measure $\nu_x$ defined on $\mathfrak{B}_Y$, such that for all $C\in \mathfrak{B}_{X\times Y}$, the function
\begin{equation}\label{meas0}
	 x \in X\setminus E_1 \longrightarrow \nu_x\left( C_x \right)  ,
	\end{equation}
where $ C_x = C\cap\left(  \left\lbrace x \right\rbrace\times Y \right)   $,	is $\mathfrak{B}_X-$measurable and 
\begin{equation}\label{int0}
	 \zeta(C)= \int_X\nu_x\left(C_x \right) \mu(dx).
	\end{equation}
	In particular,
	\[ \forall B\in \mathfrak{B}_Y, \ \nu(B)= \int_{X} \nu_x\left(B \right) \mu(dx)  ,\] 
Moreover, for all positive $\mathfrak{B}_{X\times Y}-$measurable function,  $f: X\times Y\rightarrow \mathbb{R}$, 
\begin{equation}\label{meas0f}
x \in X\setminus E_1 \longrightarrow \int_Y f(x,y)\, \nu_x\left(dy\right) 
\end{equation}
 is $\mathfrak{B}_X-$measurable   and 
 \begin{equation}\label{int0f}
	\int_{X\times Y} f(x,y)\, \zeta\left(dx,dy\right)= \int_X \left( \int_Y f(x,y)\, \nu_x\left(dy\right)  \right)  \mu(dx).
	\end{equation}
	 
\vskip 3mm
\noindent
Likewise, there is a  $\nu-$negligible set $E_2  \in \mathfrak{B}_Y$, such that   for every $\,  y\in Y\setminus E_2 $   
  there is a measure $\mu_y$  on $\mathfrak{B}_X$,  such that for all $C\in \mathfrak{B}_{X\times Y}$, the function
\[ y\in Y\setminus E_2 \longrightarrow \mu_y\left(C_y \right)  ,\]
where $ C_y = C\cap\left(  X \times\left\lbrace y \right\rbrace   \right)   $,	 is $\mathfrak{B}_{ Y}-$measurable and 
\[  \zeta(C)= \int_{X} \mu_y\left(C_y \right)  \mu(dx). \]
In particular, 
\[ \forall A\in \mathfrak{B}_X, \ \mu(A)=\int_{Y} \mu_y\left(A \right) \nu(dx)  .\]
Moreover, for all positive $\mathfrak{B}_{X\times Y}-$measurable function,  $f: X\times Y\rightarrow \mathbb{R}$, 
\[ y\in Y\setminus E_2 \rightarrow \int_{X} f(x,y)\, \mu_y\left(dx\right)  \] is  $\mathfrak{B}_{Y}-$measurable  and 
\[  \int_{X\times Y} f(x,y)\, \zeta \left(dx, dy\right) = \int_{Y} \left(  \int_{X} f(x,y)\, \mu_y \left(dx\right) \right)  \nu(dy). \]

\vskip 3mm
\noindent
As a consequence, given a $\mathfrak{B}_{X\times Y}-$measurable function,  $f: X\times Y\rightarrow \mathbb{R}$, the following affirmations are equivalent

\vskip 2mm

\begin{enumerate}
	\item $f: X\times Y\rightarrow \mathbb{R}$ is $\zeta-$integrable.
\vskip 3mm
\item $ x \in X\setminus E_1 \rightarrow \int_Y \left| f(x,y)\right| \, \nu_x\left( dy\right)   $
is $\mu-$integrable.
\vskip 3mm
\item $ y\in Y\setminus E_2 \rightarrow \int_{X}\left| f(x,y)\right| \, \mu_y\left( dx\right) $  is $\nu-$integrable.
\end{enumerate}
\vskip 2mm\noindent
And
\[ \int_{X\times Y} f(x,y)\, \zeta\left( dx, dy\right)= \int_X \left( \int_Y f(x,y)\, \nu_x\left( dy\right)  \right)  \mu\left( dx\right) . \]
\[  \qquad \qquad \qquad \qquad \qquad = \int_{Y} \left(  \int_{X} f(x,y)\, \mu_y \left( dx\right) \right)  \nu\left(  dy\right) . \]
}

\section{Proof of Theorem \ref{teo}}
\noindent
\emph{Note about the notation}: We will use $x$ and $y$ to denote generic points in $X$ and $Y$ respectively. In this way,  $B_r\left( x\right) $ automatically refers to a ball in $X$, of center $x$ and and radius $r$, while $B_r\left( y\right)$ represents a ball in $Y$ (different space, different metric). 
\vskip 2mm
\noindent
The proof of Theorem \ref{teo} will be given in several steps. 
\subsection{Definition of $ \mathit{l}_x $} 

Let  $\mathcal{Y} $ be a  dense subset of $Y$. Denote by $ \mathcal{B} $  the set of open balls $B_r\left( y\right) $ with center $y\in \mathcal{Y} $  and radius  $ r\in\mathbb{Q} $. i.e. 
\begin{equation}\label{balls}
\mathcal{B}:= \left\lbrace B_r \left( y\right) ;\, y\in\mathcal{Y},\ r\in \mathbb{Q} \right\rbrace .
\end{equation}

\noindent
Consider also the complement of the closed balls, 
\[  \mathcal{B}_{\mathrm{c}} :=\left\lbrace Y\setminus \overline{B}_r\left( y\right) ; {B}_r\left( y\right) \in\mathcal{B} \right\rbrace .  \]
\noindent
Finally, let $ \mathcal{L} $ be the set of finite unions of finite intersections of elements of $ \mathcal{B}\cup \mathcal{B}_{\mathrm{c}} $ (note that $\emptyset\in \mathcal{L}$). 
\vskip 2mm
\noindent
For each $O\in \mathcal{L}$, define the measure
\[ A\in \mathfrak{B}_X\rightarrow \mu_O \left( A\right) =\zeta\left( A\times O\right)  .\]

\noindent
Since  $\mu_O \left( A\right) \le \mu\left( A\right) $,   $\mu_O $ is absolutely continuous with respect to $ \mu $. By  Radon-Nikodym's Theorem, there is a density function $ \frac{d\mu_O}{d\mu} $, defined $\mu-$almost surely, such that $\mu_O $ is represented as an integral of this density against $\mu$. 
\vskip 2mm
\noindent
To obtain a common exceptional $\mu-$negligible set outside which $ \frac{d\mu_O}{d\mu} $ 
is well defined (by a formula) for all $ O\in \mathcal{L} $, we choose the version of $ \frac{d\mu_O}{ d\mu} $ given by the limit of the quotient of balls. To avoid talking about measurability issues, we fix once and for all a sequence $\rho_k$  decreasing to $0$. Given $ O\in \mathcal{L} $,  define
\[ \overline{\mathit{l}}_x\left( O\right) :=\limsup_{k\rightarrow \infty} \frac{\mu_O \left( B_{\rho_k}(x)\right) }{\mu \left( B_{\rho_k}(x)\right) },\]
and \[ \underline{\mathit{l}}_x\left( O\right) :=\liminf_{k\rightarrow \infty} \frac{\mu_O \left( B_{\rho_k}(x)\right) }{\mu \left( B_{\rho_k}(x)\right) },\]
\vskip 2mm
\noindent
It is well known, by a generalization of Lebesgue differentiation theorem (see for example Federer \cite{fed}, section 2.9),
 that $ \overline{\mathit{l}}_x $ and $ \underline{\mathit{l}}_x $  are versions of 
 $ \frac{d\mu_O}{d\mu} $. 
i.e. For all $  A\in \mathfrak{B}_X$,
\begin{equation}\label{int}
\mu_O \left( A\right) =\int_A  \overline{\mathit{l}}_x\left( O\right) \mu\left( dx\right)  =\int_A  \underline{\mathit{l}}_x \left( O\right) \mu\left( dx\right) .
\end{equation}
In particular,
$$ \overline{\mathit{l}}_x\left( O\right) = \underline{\mathit{l}}_x\left( O\right) ,\  \mu-\mathrm{a.s.}.$$
\noindent
Let $E_O$ be the exceptional set where the limit does not exist. i.e.
$$E_O:=\left\lbrace x\in X;\ \overline{\mathit{l}}_x\left( O\right) -\underline{\mathit{l}}_x \left( O\right) >0\right\rbrace \in \mathfrak{B}_X$$
 \noindent
Put
\[ E:=\bigcup_{O\in \mathcal{L} } E_O. \]

\noindent
Since  $\mu\left( E_O\right) =0$ for all $ O\in \mathcal{L} $, and $\mathcal{L}$ is numerable, we  have $$\mu\left( E\right) =0 .$$
 For  $O\in \mathcal{L}$, and   $x\in X\setminus E$, put 
 $$ \mathit{l}_x\left( O\right) = \overline{\mathit{l}}_x \left( O\right) =\underline{\mathit{l}}_x \left( O\right) .$$

\vskip 2pt 
\noindent
Recapitulating, for all $O\in \mathcal{L}$ and every $ A\in \mathfrak{B}_Y $, by (\ref{int}) we have
\begin{equation}\label{int1}
 \zeta\left( A\times O\right)  =\int_A  {\mathit{l}}_x\left( O\right) \mu\left( dx\right)  .
\end{equation}

\subsection{Outer measure $ \nu^*_x $} 

Changing the standpoint, we  fix   $x\in X\setminus E$ and consider the set function 
$$ O\in \mathcal{L}\rightarrow \mathit{l}_x\left( O\right) .$$
For future reference, observe that $ \mathit{l}_x $ has the following properties:
For all $ O$ and  $ \tilde{O}\in \mathcal{L} $,
\vskip 3pt
\noindent
\emph{Finite additivity}
\begin{equation}\label{add}
 \mathit{l}_x\left(  O\right) +\mathit{l}_x\left(  \tilde{O}\right) = \mathit{l}_x\left(  O\cup \tilde{O}\right) + \mathit{l}_x\left(  O\cap \tilde{O}\right) ,
\end{equation}
\vskip 2pt
\noindent
\emph{Finite subadditivity}
	\begin{equation}\label{subadd}
\mathit{l}_x\left(  O\cup \tilde{O}\right)\le	\mathit{l}_x\left(  O\right) +\mathit{l}_x\left(  \tilde{O}\right) 
	\end{equation}
\vskip 2pt
\noindent
 \emph{Monotonicity}
\begin{equation}\label{mono}	
   O\subseteq \tilde{O} \, \Rightarrow \, \mathit{l}_x\left(  O\right) \le \mathit{l}_x\left(  \tilde{O}\right)
\end{equation}
\vskip 3pt
\noindent
We need a measure on $\mathfrak{B}_Y$, capable of fulfilling the role of  $\mathit{l}_x$ in equation (\ref{int1}). Let us start by defining the outer measure
\begin{equation} \label{outer}
C\subseteq Y\rightarrow \nu^*_x\left( C\right) :=\inf\sum_{i=1}^{\infty}\mathit{l}_x\left( O_i\right) ,
\end{equation}
where the infimum is taken over all the covers $ \left\lbrace O_i \right\rbrace _{i\in \mathbb{N}} \subseteq \mathcal{L} $ of $C$. i.e. 
\[  C\subseteq \bigcup_{i=1}^\infty O_i  \, ,\ \mathrm{and}\ \forall i\in \mathbb{N}, \,  O_i\in  \mathcal{L}\]

\noindent 
It is well known that $\nu^*_x$, restricted to the set of $\nu^*_x-$measurable sets,  is a measure (denoted $\nu_x$). What we need to prove are: Firstly, that every Borel subset of $Y$ (i.e. in $\mathfrak{B}_Y$) is $\nu^*_x-$measurable and secondly that the integration property (\ref{int1}) is preserved when $\mathit{l}_x$ is replaced by $\nu_x$ (and therefore valid for any set in $\mathfrak{B}_Y$).
\vskip 2mm
\noindent
 Unfortunately, $\mathit{l}_x$ is not countably subadditive, as a  result,  $\nu^*_x$  is not  the extension of $\mathit{l}_x$, hardening  our task  a little bit. Indeed, for all $ O\in  \mathcal{L} $, we clearly have $ \nu^*_x(O)\le \mathit{l}_x(O)$, but the inverse inequality may fail, as the following example shows.
\vskip 3mm
\noindent
\emph{Example}: Let $X=Y=\left[ 0,1\right] $. Let $\mu=\nu$ be the Lebesgue measure on $ \left[ 0,1\right]  $ and $\zeta$ the normalized length on the diagonal $\left\lbrace \left(x,x \right) ;\, x\in \left[ 0,1\right] \right\rbrace $.
\vskip 3mm
\noindent
Observe that for every $ x\in \left] 0,1\right[  \, $ and $\rho_k$ small enough,
\[  \frac{\mu_{\left[ 0,x\right[} \left( B_{\rho_k}(x)\right) }{\mu \left( B_{\rho_k}(x)\right) }=\frac{1}{2}.\]

\noindent
So, $ \mathit{l}_x  \left( \left[ 0,x\right[\right) =   1/2 $, while $ \nu_x^*  \left( \left[ 0,x\right[\right) =   0 $. In fact, we can cover $ \left[ 0,x\right[ $ with a sequence of intervals $ \left[ 0,x_n\right[ $, where $x_n \in \mathbb{Q}$ increases to $x$. Each one of the intervals $ \left[ 0,x_n\right[ $ verify
\[  \frac{\mu_{\left[ 0,x_n\right[} \left( B_{\rho_k}(x)\right) }{\mu \left( B_{\rho_k}(x)\right) }=0, \]
for all  $\rho_k$ small enough. Therefore, for all $n\in \mathbb{N}$,
$ \mathit{l}_x  \left( \left[ 0,x_n\right[\right) =0 $ and
\[ \nu_x^*  \left( \left[ 0,x\right[\right) \le \sum_{i=1}^\infty  \mathit{l}_x  \left( \left[ 0,x_n\right[\right) =   0 .\]
\vskip 3mm
\noindent 
\subsection{Borel subsets of $Y$ are $\nu^*_x-$measurable} Let's prove first that any open ball $B_r(y)\in\mathcal{L}$ is $\nu^*_x-$measurable. To this end, fix $C\subseteq Y$ and a cover  $ \left\lbrace O_i\right\rbrace_{i\in\mathbb{N}}  \subseteq \mathcal{L}$ of $C$. We must  show that
\begin{equation}\label{measurable}
\nu^*_x \left( C \cap B_r(y) \right) + \nu^*_x \left( C\setminus B_r(y)\right) \le \sum_{i=1}^\infty \mathit{l}_x\left( O_i\right).
\end{equation}
\noindent
Since $ O_i\cap B_r(y)  \in \mathcal{L}$ and  $ \left\lbrace O_i\cap B_r(y)  \right\rbrace_{i\in\mathbb{N}} $  is a covering of $ C\cap B_r(y) $, 
\begin{equation}\label{uno}
 \nu^*_x\left( C\cap B_r(y)\right) \le \sum_{i=1}^\infty  \mathit{l}_x\left( O_i\cap B_r(y)  \right)    .
\end{equation}
\noindent
Now, let $\alpha_{i}<1$, $\alpha_{i}\in\mathbb{Q}$. Then,  $O_i\setminus \overline{B}_{ \alpha_{i} r}(y)\in \mathcal{L}$ and  $ \left\lbrace O_i\setminus \overline{B}_{ \alpha_{k_i} r}(y) \right\rbrace_{i\in\mathbb{N}}  $ is a covering of $C\setminus B_r(y) $. So,
  \begin{equation}\label{dos}
 \nu^*_x\left( C\setminus B_r(y)\right) \le \sum_{i=1}^\infty  \mathit{l}_x\left(   O_i\setminus \overline{B}_{ \alpha_i r}(y)  \right) .
 \end{equation}
 By (\ref{add}) and (\ref{mono}), we have
\begin{equation}\label{tres} 
 \mathit{l}_x\left( O_i\cap B_r(y)  \right)  +   \mathit{l}_x\left(   O_i\setminus \overline{B}_{ \alpha r}(y)  \right) \le  \mathit{l}_x\left( O_i \right)  +   \mathit{l}_x \left( B_r(y) \setminus \overline{B}_{ \alpha_i r}(y)  \right).  
 \end{equation}
\vskip 2mm
\noindent
Adding (\ref{uno}) and (\ref{dos}), and using (\ref{tres}), we obtain
\[\nu^*_x\left( C\cap B_r(y)\right) +\nu^*_x\left( C\setminus B_r(y)\right) \]
\[\le \sum_{i=1}^\infty  \mathit{l}_x\left( O_i \right) +\sum_{i=1}^\infty   \mathit{l}_x\left( B_r(y)\setminus \overline{B}_{ \alpha_i r}(y)  \right) .\]
\vskip 2mm
\noindent
The result follows if we can make 
the second sum as small as we want. 
Unfortunately, for a fixed $x\in E$, we might fail to do so, even though  
\begin{equation}\label{empty}
\bigcap_{i=1}^\infty B_r(y) \setminus \overline{B}_{ \alpha_i r}(y) =\emptyset ,
\end{equation} 
for any sequence  $\alpha_i\rightarrow 1$. 
In fact,  we can not switch the limits in
\[ \lim_{i\rightarrow \infty}\mathit{l}_x\left( B_r(y)\setminus \overline{B}_{ \alpha_i r}(y)  \right) = \lim_{i\rightarrow \infty}\lim_{k\rightarrow \infty} \frac{\zeta \left( B_{\rho_k}\left( x\right) \times \left(   B_r(y)\setminus \overline{B}_{ \alpha_i r}(y)  \right) \right) }{\mu\left(  B_{\rho_k}\left( x\right) \right)} .\]
So, we need to look back at what happens for $x$ variable and check whether we can solve the problem by  throwing  away a few more points (meaning to enlarge $E$).
\vskip 2mm
\noindent
By (\ref{int1}) and ($\ref{empty}$), and $\gamma < 1$, 
\[ \int_X \mathit{l}_x\left( B_r(y) \setminus \overline{B}_{ \gamma r}(y)\right)  \mu\left( dx\right) =\nu\left(  B_r(y) \setminus \overline{B}_{ \gamma r}(y)\right) \longrightarrow_{\gamma \nearrow 1} 0. \]
Therefore, 
\begin{equation}\label{null-a.e.}
\mathit{l}_x\left( B_r(y) \setminus \overline{B}_{ \gamma r}(y)\right)  \longrightarrow_{\gamma \nearrow 1} 0,\ \mu-\mathrm{a.s.}  
\end{equation}
\vskip 2mm
\noindent
Now, fix an increasing sequence  $ \left\lbrace \gamma_j\right\rbrace_{j\in\mathbb{N}}  \subseteq  \mathbb{Q}$, $ \gamma_j \rightarrow 1$. Define, for all $y\in\mathcal{Y}$ and $r\in\mathbb{Q}$,
\[ {E}_{r,y}:=\left\lbrace x\in X\setminus E ;\ \liminf_{j\rightarrow \infty}\mathit{l}_x\left( B_r(y) \setminus \overline{B}_{  \gamma_j r}(y)\right)  > 0 \right\rbrace  .\]
By (\ref{null-a.e.}), $\mu\left( {E}_{r,y}\right) =0 $. Since
\[   {E}_1:= E\cup\bigcup_{r\in\mathbb{Q} ,\, y\in\mathcal{Y} }{E}_{r,y} \]
is a countable union of sets of $\mu-$measure 0, we  have $$\mu\left( {E}_1\right) =0 . $$
For all  $x\in X\setminus E_1  $, and every $\epsilon>0$, we can choose a subsequence $ \left\lbrace \alpha_i\right\rbrace_{i\in\mathbb{N}}  $ of $ \left\lbrace \gamma_j\right\rbrace_{k\in\mathbb{N}}  $, such that
\[ \sum_{i=1}^\infty   \mathit{l}_x\left( B_r(y)\setminus \overline{B}_{ \alpha_i r}(y)  \right) < \epsilon .\]
This completes the proof of (\ref{measurable}). 
\vskip 2mm
\noindent
Consequently, for all  $x \in X\setminus   E_1$, $\nu_x$ is a measure defined at least in the $\sigma-$field generated by $\mathcal{B}$, i.e.  $\mathfrak{B}_Y$.

\subsection{Measurability and integrability  for compact sets}
Our task now is to prove that given $B\in \mathfrak{B}_Y$, the function
\begin{equation}\label{meas}
 x\in X\setminus  E_1 \longrightarrow \nu_x\left( B\right)  
\end{equation} 
is measurable and, for all $A\in \mathfrak{B}_X$
\begin{equation}\label{int2}
\zeta\left( A\times B\right) =\int_A  \nu_x\left( B\right) \mu\left( dx\right)   .
\end{equation}
\vskip 2pt
\noindent
Let's consider first a finite intersection of closed, compact balls  
$$ \overline{B}=\overline{B}_{r_1}(y_1)
\cap\cdots\cap \overline{B}_{r_n}(y_n) ,$$
with $y_1,\cdots,y_n\in\mathcal{Y}$  and $ r_1,\cdots,r_n\in\mathbb{Q} $.
 By the measurability of $\mathit{l}_x$ and (\ref{int1}), the properties (\ref{meas}) and (\ref{int2}) are proven at once if we show  
\begin{equation}\label{l_x}
\nu_x \left( \overline{B} \right) =1-\mathit{l}_x\left( Y\setminus \overline{B}\right) ,\ \mu-\mathrm{a.s.}  
\end{equation}
(We use $\mathit{l}_x\left( Y\setminus \overline{B}\right) $ just because $\mathit{l}_x $ is not defined for $ \overline{B}$.)
\vskip 2pt
\noindent
Let $O_1,\, O_2,  \cdots, O_m\subseteq \mathcal{L}$ be a covering of 
$ \overline{B} $. We can assume the covering is finite, since $\overline{B}$ is compact and the sets in $\mathcal{L}$ are open.
\vskip 3pt
\noindent
Clearly, for all $ x\in X\setminus E $,
\[ \lim_{k\rightarrow \infty} \frac{\zeta\left( {B}_{\rho_k}(x)\times \overline{B}\right) }{\mu\left( {B}_{\rho_k}(x)\right) }=\lim_{k\rightarrow \infty} \left( 1-\frac{\zeta\left( {B}_{\rho_k}(x)\times \left( Y\setminus \overline{B} \right)\right)  }{\mu\left( {B}_{\rho_k} (x)\right) }\right) =1-\mathit{l}_x\left( Y\setminus \overline{B}\right) . \] 
Since the covering is finite, by  (\ref{subadd}),
 \[ \lim_{k\rightarrow \infty}  \frac{\zeta\left( {B}_{\rho_k}(x)\times \overline{B}\right) }{\mu\left( {B}_{\rho_k}(x)\right) }\le \mathit{l}_x\left( \bigcup_{i=1}^{n} O_i\right) \le \mathit{l}_x\left( O_1\right) + \cdots +   \mathit{l}_x\left( O_m\right).  \] 
Then,
\[ 1-\mathit{l}_x\left( Y\setminus \overline{B}\right) \le  \nu_x \left( \overline{B} \right) .\] 
\vskip 2pt
\noindent
On the other hand, given $ \eta >1 $, $ \eta\in\mathbb{Q} $, and denoting
\[ B_\eta={B}_{\eta r_1}(y_1)
\cap\cdots\cap {B}_{\eta r_n}(y_n) , \]
we have 
\[ \nu_x \left( \overline{B}\right) \le \mathit{l}_x\left(  {B}_\eta \right) .\] 
Consequently, taking any sequence $\eta_k\searrow 1$, $ \eta_k \in\mathbb{Q} $, we have
\[ 1-\mathit{l}_x\left( Y\setminus \overline{B}\right) \le  \nu_x \left( \overline{B}\right)  \le \lim_{k\rightarrow\infty}\mathit{l}_x\left(  {B}_{\eta_k}\right) .\] 
\vskip 2pt
\noindent
Since the functions at the left and at the  right of the above inequalities are $\mathfrak{B}_{X}-$measurable, and equal between them $\mu-$almost-surely, using (\ref{int1}), we have proved (\ref{l_x})  and, a fortiori,  (\ref{meas})  and (\ref{int2}), at least for a finite intersection of compact balls.

\subsection{Measurability and integrability  for Borel sets}
 Let $\mathcal{M}$ be the collection  of sets  $B\in \mathfrak{B}_Y$  actually verifying (\ref{meas})  and  (\ref{int2}). 
 \vskip 3pt
 \noindent
 Clearly, $\emptyset\in \mathcal{M} $ and $ Y\setminus B\in \mathcal{M} $, whenever  $B \in \mathcal{M}  $. 
\vskip 3pt
\noindent
Now, take a disjoint sequence   $B_1,B_2,\, \cdots \, \in \mathcal{M}$  ($B_i\cap B_j= \emptyset $,   $ i\neq j $).  
Since, for all $x\in X\setminus  E_1  $, $\nu_x$ is a measure on $\mathfrak{B}_Y$,
\begin{equation}\label{suma}
 \nu_x\left( \bigcup_{i=1}^{\infty} B_i\right) = \sum_{i=1}^{\infty}\nu_x\left(  B_i\right) . 
\end{equation}
By  (\ref{meas}), $ x\rightarrow  \nu_x\left( \bigcup_{i=1}^{\infty} B_i\right)  $ is $\mathfrak{B}_{X}-$measurable, being a countable sum of $\mathfrak{B}_{X}-$measurable functions. 
\vskip 2pt
\noindent
By (\ref{suma}), the monotone convergence theorem, and  (\ref{int2}) (remember $B_i \in \mathcal{M}$, for all  $i\in \mathbb{N}$), given $ A\in\mathfrak{B}_X $,
\[ \int_A  \nu_x\left( \bigcup_{i=1}^{\infty} B_i\right) \mu\left( dx\right) = \int_A  \sum_{i=1}^{\infty}\nu_x\left(  B_i\right) \mu\left( dx\right)   \]
\[ \ \quad  \qquad  \qquad  \qquad  \qquad = \sum_{i=1}^{\infty} \int_A  \nu_x\left(  B_i\right)  \mu\left( dx\right)  \]
\[     \qquad   \quad  \qquad  \qquad = \sum_{i=1}^{\infty}   \zeta\left( A\times B_i\right) \]
\[ \qquad \qquad  \qquad  \qquad = \zeta \left( \bigcup_{i=1}^{\infty} A\times B_i\right) . \]
\[ \qquad \qquad  \qquad  \qquad = \zeta \left( A\times\bigcup_{i=1}^{\infty}  B_i\right) . \]
 Then, 
 \[ \bigcup_{i=1}^{\infty} B_i\in \mathcal{M}. \]
Since $\mathcal{M}$ contains all the finite intersections of compact balls with center in $\mathcal{Y}$ and rational radius, by the $\pi-\lambda$ theorem (see \cite{Bi}, page 36),  
 $$ \mathcal{M}=\mathfrak{B}_Y. $$

\subsection{Proof of the theorem} We have proven so far (\ref{meas0}) and (\ref{int0})   for sets of the form $C=A\times B$, with $ A\in \mathfrak{B}_X $  and $ B\in \mathfrak{B}_Y $. 
\vskip 2mm
\noindent
Using a similar argument as before, let $\tilde{\mathcal{M}}$ denote the collection  of sets  $C\in \mathfrak{B}_{X\times Y} $  verifying (\ref{meas0})  and  (\ref{int0}). We readily see that $\emptyset\in \mathcal{M} $. 

\vskip 2mm
\noindent
Now, let $C\in\tilde{\mathcal{M}} $. Since $\left( \left( X\times Y\right) \setminus  C\right)_x = Y  \setminus  C_x $, the function
\[ x\rightarrow \nu_x\left(\left( \left( X\times Y\right) \setminus  C\right) _x \right)  = 1- \nu_x\left( C_x\right) \]
 is $ \mathfrak{B}_{X}-$measurable, and 
\[ \int_X \nu_x\left(\left( \left( X\times Y\right) \setminus  C\right)_x \right)  \mu\left( dx\right)  = 1- \int_X \nu_x\left(C_x \right) \mu\left( dx \right)   \]
\[ \ \qquad  \qquad  \qquad  \qquad = 1-\zeta\left(  C\right)\]
\[\  \qquad  \qquad  \qquad  \qquad  \qquad  \ \quad = \zeta\left( \left( X\times Y\right) \setminus  C\right). \]
\noindent
Then, for all $C\in\tilde{\mathcal{M}} $, we have $\left( X\times Y\right) \setminus  C \in\tilde{\mathcal{M}} $.
\vskip 3pt
\noindent
Finally,  let  $C_1,C_2,\, \cdots \, \in \tilde{\mathcal{M}} $ a sequence of disjoint sets ($C_i\cap C_j= \emptyset $, for all $ i\neq j $).  
Since 
$\nu_x$ is a measure on $\mathfrak{B}_Y$,
\begin{equation}\label{suma1}
\nu_x\left(  \left( \bigcup_{i=1}^{\infty} C_i\right)_x\,\right) = \nu_x\left( \bigcup_{i=1}^{\infty} \left( C_i\right) _x\right) = \sum_{i=1}^{\infty}\nu_x\left(  \left( C_i\right)_x \right)  . 
\end{equation}
Then, $  x\rightarrow  \nu_x\left(  \left( \bigcup_{i=1}^{\infty} C_i\right)_x \,\right)  $ is $\mathfrak{B}_{X}-$measurable, being a sum of $\mathfrak{B}_{X}-$measurable functions, since $C_i \in \tilde{\mathcal{M}}$,  for all $i\in \mathbb{N}$.   

\vskip 3pt\noindent
By (\ref{suma1}), the monotone convergence theorem, and  (\ref{int0})  ($C_i \in \tilde{\mathcal{M}}$), 
\[ \int_X  \nu_x\left(  \left( \bigcup_{i=1}^{\infty} C_i\right)_x\,\right) \mu\left( dx \right) = \int_X  \sum_{i=1}^{\infty}\nu_x\left(  \left( C_i\right) _x\right) \mu\left( dx \right) \ \quad   \quad   \]
\[ \ \quad  \qquad  \qquad  \qquad  \qquad = \sum_{i=1}^{\infty} \int_X  \nu_x\left(  \left( C_i\right) _x\right)  \mu\left( dx \right) \]
\[  \qquad  \qquad = \sum_{i=1}^{\infty}   \zeta\left(  C_i\right) \]
\[ \  \qquad  \qquad = \zeta \left( \bigcup_{i=1}^{\infty}  C_i\right)  \]
Then, 
\[ \bigcup_{i=1}^{\infty} C_i\in \tilde{\mathcal{M}}. \]
\vskip 2mm
\noindent
Since 
\[ \left\lbrace A\times B;\, A\in \mathfrak{B}_X ,\, B\in \mathfrak{B}_Y \right\rbrace  \subseteq \tilde{\mathcal{M}} ,\]
by the $\pi-\lambda$ theorem,
\[ \tilde{\mathcal{M}} = {\mathfrak{B}_{X\times Y}}  .\]
\vskip 3mm
\noindent
The remaining of the proof is standard. Assume $f :{X\times Y}\rightarrow \mathbb{R}$ is a  positive, $\mathfrak{B}_{X\times Y}-$measurable functions. Then, $f$ can be approached by an increasing sequence of simple, $\mathfrak{B}_{X\times Y}-$measurable functions  (linear combinations of characteristic functions)
\[ f_k(x,y)=\sum_{i=1}^{n_k} \lambda_{k,i}C_{k,i},\]
where $\lambda_{k,i}\in\mathbb{R}$ and  $C_{k,i} \in \mathfrak{B}_{X \times Y}$ ($i=1,\cdots, n_k $, $ k\in  \mathbb{N}$).
\vskip 3mm
\noindent
By linearity, for each $ f_k$, properties (\ref{meas0f}) and (\ref{int0f}) are a direct consequence of (\ref{meas0}) and (\ref{int0}). Passing to the limit as $k\rightarrow \infty$, equations (\ref{meas0f}) and (\ref{int0f}) are conserved, therefore valid  for any positive function like $f$.
\vskip 3mm
\noindent
To establish the integrability equivalence, we apply the preceding result to the positive and negative parts of the given function. In this way,  the proof of Theorem \ref{teo} is complete.

\section{An aplication to Optimal Transport}  \label{OT}
\noindent
In this section we use Theorem \ref{teo} to show that a pair of competitive price functions, whose integral with respect to some transference plan equals the transport cost, are  conjugate to each other almost surely, complementing  Kantorovich's duality theorem on the nature of a pairs of competitive prices maximizing the profit. See Villani's book \cite{vill}, page 70, for a very detailed discussion on  Kantorovich's theorem, in particular  
 Theorem 5.1, part \textit{(ii)}, item \textit{(d)}.  For this result, we do not need lower semicontinuity  of the cost and other assumptions used to prove Kantorovich's theorem, so we state the following  lemma in its simplest form, using the notations in  Theorem \ref{teo}, by the way.

{\lem\label{ot} Let $ c : X \times Y \rightarrow \mathbb{R} \cup \left\lbrace +\infty\right\rbrace $ be a ${\mathfrak{B}_{X\times Y}}-$measurable, 
	cost function, and
	$ \psi :X\rightarrow \mathbb{R}\cup \left\lbrace +\infty\right\rbrace  \ \mathrm{and}\ \phi :Y\rightarrow \mathbb{R}\cup \left\lbrace -\infty\right\rbrace $
	be a pair of  competitive prices, i.e. 
	\[ \forall \, \left( x,y\right) \in X\times Y, \ \phi \left( y\right) -\psi\left( x\right) \le c\left( x,y\right) . \]
Let  $\pi$  be a probability measure on ${\mathfrak{B}_{X\times Y}}$, with  marginal measures $\mu$ on ${\mathfrak{B}_X}$ and $\nu$ on ${\mathfrak{B}_Y}$. (i.e.  $\pi$   is a transference plan between $\mu$ and $\nu$.)
Assume  that $ \phi -\psi$ and $ c $ are $\pi-$integrable and 
\[  \phi \left( y\right) -\psi\left( x\right) = c\left( x,y\right) , \ \pi-\mathrm{a.s.} . \]
Then 
\begin{equation}\label{intes01}
\psi\left( x\right) = \sup_{y\in Y}\left(\phi \left( y\right)  - c\left( x,y\right) \right)  , \ \mu-\mathrm{a.s.} 
\end{equation}
and 
\begin{equation}\label{intes02}
 \phi \left( y\right) = \inf_{x\in X}\left( \psi \left( y\right)  + c\left( x,y\right)\right) , \ \nu-\mathrm{a.s.} 
 \end{equation}
}
\vskip 2mm
\noindent
\emph{Proof}: Since $ \psi+ c-\phi = 0 $, $\pi-$a.s., 
\[\int_{X\times Y} \left( \psi(x)+ c(x,y)-\phi(y) \right)  \pi\left(dx,dy\right)= 0. \]
Using the decomposition  given by Theorem \ref{teo}, equation (\ref{int0f}), 
\begin{equation}\label{nullint}
\int_X \left( \int_Y \left( \psi(x)+ c(x,y)-\phi(y) \right)  \nu_x\left(dy\right)  \right)  \mu(dx)=0. 
\end{equation}
Since $ \psi+ c-\phi \ge 0 $, 
\begin{equation}\label{sup}
\psi(x)\ge \sup_{y\in Y}  \left(\phi(y) - c(x,y) \right) 
\end{equation}
and, for all $x$ where it is defined, the function
\[ x\longrightarrow \int_Y \left( \psi(x)+ c(x,y)-\phi(y) \right)  \nu_x\left(dy\right)  
 \]
is nonnegative. By (\ref{nullint}), 
\[  \int_Y \left( \psi(x)+ c(x,y)-\phi(y) \right)  \nu_x\left(dy\right)  = 0,\ \mu-\mathrm{a.s.} \]
Then,
 \begin{equation}\label{intas}
\psi(x)= \int_Y \left(\phi(y) - c(x,y)\right)  \nu_x\left(dy\right)\le \sup_{y\in Y}  \left(\phi(y) - c(x,y) \right) ,\ \mu-\mathrm{a.s.}
\end{equation}
Combining (\ref{sup}) and (\ref{intas}), we obtain (\ref{intes01}). Equation (\ref{intes02}) is validated in a similar way.

\bibliographystyle{amsplain}
\bibliography{fubini01}

\providecommand{\bysame}{\leavevmode\hbox to3em{\hrulefill}\thinspace}
\providecommand{\MR}{\relax\ifhmode\unskip\space\fi MR }
\providecommand{\MRhref}[2]{%
  \href{http://www.ams.org/mathscinet-getitem?mr=#1}{#2}
}
\providecommand{\href}[2]{#2}
\begin{thebibliography}{1}

\bibitem{Bi}
Patrick Billigsley, \emph{Probability and measure}, Second edition, John Willey
  \& Sons, 1986.

\bibitem{fed}
Herbert Federer, \emph{Geometric measure theory}, Springer-Verlag, Berlin -
  Heidelberg - NewYork, 1969.

\bibitem{sasaki}
Shigeo Sasaki, \emph{On differential geometry of tangent bundles of riemannian
  manifolds}, Tohoku Math. J. (2) \textbf{10, Number 3} (1958), 338--354.

\bibitem{vill}
C\'edric Villani, \emph{Optimal transport, old and new}, Springer, Berlin,
  2008.

\end{thebibliography}

\end{document}